\documentclass [12pt,a4paper,reqno]{amsart}
\usepackage{amsmath,amsthm}
\usepackage{amscd}
\usepackage{amsfonts}
\input amssymb.sty
\chardef\bslash=`\\ 

\newtheorem{thm}{Theorem}
\newtheorem*{thm*}{Theorem}
\newtheorem*{conjecture*}{Conjecture}

\newtheorem{cor}{Corollary}

\newtheorem{lem}{Lemma}
\newtheorem{prop}{Proposition}

\theoremstyle{definition}
\newtheorem{defn}{Definition}
\newtheorem*{remark*}{Remarks}
\newtheorem*{examples*}{Examples}
\newtheorem*{defn*}{Definition}
\newtheorem*{cor*}{Corollary}

\newcommand{\thmref}[1]{Theorem~\ref{#1}}
\newcommand{\secref}[1]{Section~\ref{#1}}

\newcommand{\lemref}[1]{Lemma~\ref{#1}}
\newcommand{\corref}[1]{Corollary~\ref{#1}}
\newcommand{\propref}[1]{Proposition~\ref{#1}}
\newcommand{\defnref}[1]{Definition~\ref{#1}}
\date{}

\begin{document}
\title{Automorphisms of the semigroup of all endomorphisms of free algebras}
\author[Grigori Zhitomirski]{G. Zhitomirski\\
 Department of Mathematics\\
 Bar-Ilan University, 52900 Ramat Gan, Israel}
\thanks {This research is partially supported by THE ISRAEL SCIENCE FOUNDATION
founded by The Israel Academy of Sciences and Humanities - Center of Excellence Program.}

\maketitle
\begin{abstract} Last years a number of papers were devoted  to describing automorphisms of semigroups of
endomorphisms of free finitely generated universal algebras of some varieties: groups, semigroups, associative
commutative algebras, inverse semigroups, modules, Lie algebras and some others (see references). All these
researches were inspired by the questions  prof. B. Plotkin set in connection with so called universal algebraic
geometry \cite{Plot7L}, \cite{PlotAlgGeom}. The aim of this paper is to suggest a method of describing $Aut
\,End (F)$ for a free algebra $F$ of an arbitrary variety of universal algebras. This method allows to obtain
easily all known results as well as new ones.
\end{abstract}
\section{introduction}
The aim of this paper is to suggest a method of describing $Aut \,End (F)$ where $F$ is a  free algebra of an
arbitrary variety of universal algebras. The interest to this problem was inspired in most cases by questions
which were set by B. Plotkin in connection with universal algebraic geometry developed in his papers (see for
example \cite{Plot7L}, \cite{PlotAlgGeom}).

To this moment, there are many papers devoted to the mentioned problem for concrete varieties, namely, for
varieties of semigroups \cite{MashSch}, groups \cite{Form}, inverse semigroups \cite{MashSchZh}, linear algebras
\cite{Berz}, \cite{BerzPlot},\cite{KatsLipPlot}, \cite{LipPlot}, \cite{MashPlotPlot}, \cite{MashAssAlg}  and
some others. Since the considered varieties are different, the applied methods are different too, however the
obtained descriptions are very similar. This situation requires an explanation, but as far as the author knows,
there are no attempts to observe this problem from a general point of view.

In the present paper, an approach is offered that allows to solve the problem in general and to show what to do
in every concrete case. First of all, automorphisms of  $End (F)$ are described for arbitrary finitely generated
free universal algebras $F$ without any restriction on the varieties. This description can be simplified in many
important cases, in these cases all considered automorphisms are quasi-inner. An automorphism $\Phi$ is called
{\it quasi-inner} if it acts as follows: $\Phi (\nu )=s\circ \nu \circ s\sp {-1}$, where $s$ is a permutation of
$F$ which may not be an automorphism of $F$. $\Phi$ is called {\it inner} if this permutation $s$ can be chosen
as an automorphism of $F$. Necessary and sufficient conditions are given for $\Phi $ to be quasi-inner.

Roughly speaking, the approach consists of two steps. The first step is to verify that an automorphism $\Phi$ of
$End (F)$ is quasi-inner. In order to do it, some conditions are given that are necessary and sufficient for
this fact. In many cases, some of these conditions can be more or less easily checked. The second step is more
simple. We consider a quasi-inner automorphism as an inner one with respect to more wide set of morphisms than
standard homomorphisms for this variety form and then try to describe that set of morphisms. It turns out that
every quasi-inner automorphism determines by means of the mentioned above permutation $s$ an isomorphic
derivative structure on $F$. Those structures can be found and hence $s$ can be described.

Using this method, we obtain in more simple way all known results for varieties mentioned above as well as new
ones. We consider finitely generated free algebras only, but it is not a necessary restriction.

\subsection{Acknowledgements} The author is happy to thank  B.Plotkin, E.Plotkin, G.Mashevitzky, R. Lipyanski and
A. Tsurkov for the stimulating discussions of the results.

\section{main ideas and results}
Let $F=F(x\sb 1 ,\dots ,x\sb n)$ be a free algebra in a variety $V$ over the set $X=\{x\sb 1 ,\dots ,x\sb n \}$.
We consider $Aut\,End(F)$ and try to describe this group. Let $\Phi \in Aut\,End(F)$. First of all, we describe
how $\Phi $ acts of $End(F)$ (\thmref{AutoAct}, \thmref{CondPerm}). It turns out that its action is produced by
a special permutation on the direct power of $F: \; F \sp n$. In some cases, it is produced by an permutation of
$F$, that is, $\Phi$ is quasi-inner (\corref{DirProd}). In what follows, we concentrate our attention on
quasi-inner automorphisms.

The first main idea is simple. Let $\nu$ be an endomorphism of $F$. Denote by $G\sb {\nu}$ the invariant group
of $\nu$, namely, $G\sb {\nu}=\{\theta \in Aut (F)\; \vert \nu \circ \theta =\nu \}$. Conversely, for every
automorphism $\theta $ of $F$ we have a left ideal $I\sb {\theta }$ consisting of all endomorphisms $\nu$ such
that $\nu \circ \theta =\nu$. Extending this notions up to subsets of $F$ we obtain a Galois corresponding
between subgroups and left ideals of $End (F)$. Every automorphism $\Phi $ of $End (F)$ preserves this
corresponding, particularly, if $\Phi (\nu )=\mu $ then $\Phi$ maps isomorphically the group $G\sb {\nu}$ onto
group $G\sb {\mu}$, and if $\Phi (\gamma) =\theta$ for an automorphism $\gamma$ of $F$ then $\Phi$ maps
isomorphically the left ideal $I\sb {\gamma }$ on the left ideal $I\sb {\theta }$.

Let $P$ be the set of all automorphisms of $F$  that preserve $X$, i.e., they act on $X$ as permutations.
Clearly, $P$ is a subgroup of the group $Aut (F)$. We select this subgroup as very important one. Using the idea
mentioned above, we obtain conditions equivalent to the condition that an automorphism is quasi-inner
(\thmref{potinner}). In \secref{examples}, we use those conditions for some varieties. To verify them for other
cases, we introduce so called basis matrices and study properties of such matrices. The facts declared above are
the contents of \secref{secQI} and \secref{examples}.

The second idea is a development of the approach that was successfully used in \cite{PlotZhCat},
\cite{PlotZhAlg} for categories and was never used in the considered situation. It lies in studying of
polynomial operations which determine on $F$ a structure isomorphic to the original one. This idea is realized
in \secref{DerOper}. If an automorphism $\Phi$ is quasi-inner, it is produced by a permutation $s$ of $F$. It is
shown that we can restrict our considerations to the case when this permutation preserves all basis elements.
Being a permutation, $s$ translates every operation to an isomorphic one. On the other hand, the image of the
atomic term corresponding to the chosen operation is also a term and hence determines a new operation (a
derivative or a polynomial operation). It turns out that these two new operations corresponding to the chosen
one coincide \thmref{IsomDeriv}.

In the last \secref{Appl}, we consider all known facts from the presented point of view and obtain some new
ones.

\section{Quasi-inner automorphisms}\label{secQI}
We begin with an usual generalization of the notion of inner automorphism.
\begin{defn} An automorphism $\Phi$ of $End(F)$ is called {\it quasi-inner} if there exists a bijection $s$ of the set
$F$ such that $\Phi (\nu)=s\circ \nu \circ s\sp {-1}$ for every endomorphism $\nu$ of $F$. If $s$ turns out to
be an automorphism of $F$, $\Phi$ is an {\it inner} automorphism of $End(F)$.
\end{defn}

In most interesting cases, we deal with  quasi-inner automorphisms, but we have no reasons to think that it is
the general case.\footnote {R. Lipyanski has informed me that he has an example of variety with free algebras
$F$ admitting non quasi-inner automorphisms of $End(F)$.} Therefore, we consider the general case when we do not
know whether an automorphism $\Phi$ is quasi-inner. Take the direct product $F\sp n =F\times F\times \dots
\times F$ of $n$ copies of $F$ and denote its elements as vectors $\vec{u} =(u\sb 1 ,\dots ,u\sb n )$. Let
$\Delta $ be the diagonal of $F\sp n$. Of course, $\Delta $ is isomorphic to $F$ and the following natural map
$\delta $ is an isomorphism of $F$ onto $\Delta $ : $\delta (u)=(u ,\dots ,u )$.

To every endomorphism $\nu$ of $F$, an endomorphism $\nu \sp {\times n}$ is assigned as follows: $ \nu \sp
{\times n} =\nu \times \dots \times \nu \;$. Denote by $\alpha \sb {(u\sb 1 , \dots ,u\sb n)}$ an endomorphism
defined by its values on basis elements: $\alpha \sb {(u\sb 1 , \dots ,u\sb n)} (x\sb i )=u\sb i$ for $i=1,
\dots , n$. It is easy to see that we have
\begin{equation}\label{alphas}
\nu \circ \ \alpha \sb {(u\sb 1 , \dots ,u\sb n)}=\alpha \sb {(\nu (u\sb 1) , \dots ,\nu (u\sb n))}
\end{equation}
for every endomorphism $\nu $ of $F$

Let $\Phi $ be an automorphism of $End (F)$. Define a permutation $\varphi $ of $F\sp n$ corresponding to
$\Phi$:
\begin{equation}\label{genmap}
\varphi (\vec{u})=\vec{v} \Leftrightarrow \Phi (\alpha \sb {\vec{u}})=\alpha \sb {\vec{v}}.
\end{equation}
\begin{lem}\label{DeltaCondition}
If $\varphi (\vec{u})\in \Delta$ then $\varphi (\nu \sp {\times n} (\vec{u}))\in \Delta$ for every endomorphism
$\nu $ of $F$.
\end{lem}
\begin{proof} Let $\varphi (u\sb 1 , \dots ,u\sb n)=(v,\dots ,v)$ and
$\varphi (\nu (u\sb 1) , \dots ,\nu (u\sb n))=(w\sb 1 , \dots ,w\sb n)$. The last one means that $\Phi (\nu
\circ \alpha \sb {(u\sb 1 , \dots ,u\sb n)})=\alpha \sb {(w\sb 1 , \dots ,w\sb n)}$. Hence $\Phi (\nu ) \circ
\alpha \sb {(v,\dots ,v)}=\alpha \sb {(w\sb 1 , \dots ,w\sb n)}$, that implies $\Phi (\nu )(v)=w\sb i$ for all
$i=1,\dots ,n$. Thus $w\sb 1 =\dots =w\sb n$.
\end{proof}
\begin{thm}\label{AutoAct} Every automorphism $\Phi $ of $F$ can be presented in the following form:
$$\Phi (\nu )=\delta \sp {-1} \circ \varphi \circ \nu \sp {\times n} \circ \varphi \sp {-1} \circ \delta $$
for every endomorphism $\nu $ of $F$.
\end{thm}
\begin{proof} Firstly, this formula is correct because of Lemma above. Apply the endomorphisms $\Phi $ to the
both sides of the equality~(\ref{alphas}) and obtain $\Phi (\nu ) \circ \alpha \sb {\varphi (\vec{u})}=\alpha
\sb {\varphi (\nu \sp {\times n} (\vec{u}))}$. It means that $(\Phi (\nu )) \sp {\times n} =\varphi \circ \nu
\sp {\times n} \circ \varphi \sp {-1}$. This implies the required formula.
\end{proof}
Of course, the formula above does not determine an automorphism of $End(F)$ for every permutation $\varphi $ of
$F\sp n$. This permutation must satisfy some conditions which we formulate below. The algebra $F\sp n$ is a
direct power of the free algebra $F$. Any permutation $\varphi $ of $F\sp n$ determines a new algebraic
structure $(F\sp n )\sp *$ on the set $F\sp n$ such that $\varphi : (F\sp n )\sp * \to F\sp n$ is an an
isomorphism.  The obtained algebra $(F\sp n )\sp *$  is also a direct product of $n -$generated free algebras
but not necessary a direct power of an algebra. Using this new structure we can formulate a necessary and
sufficient condition for $\varphi$ to determine an automorphism $\Phi $ of $F$.
\begin{thm}\label{CondPerm}
A permutation $\varphi$ of $F\sp n$ determines an automorphism $\Phi $ of $End(F)$ if and only if, for every
endomorphism $\nu $ of $F$, $\varphi$ satisfies the condition formulated in \lemref{DeltaCondition} and the map
$\nu \sp {\times n}$ is an endomorphism of the algebra $(F\sp n )\sp *$.
\end{thm}
\begin{proof} The necessity of the condition follows from \lemref{DeltaCondition} and the expression: $\nu \sp n = \varphi \sp {-1}
\circ (\Phi (\nu ))\sp n \circ \varphi$. Let us proof that it is sufficient.

Assume that every endomorphism $\nu $ of the algebra $F$ provides that the map $\nu \sp {\times n}$ is an
endomorphisms of the algebra $(F\sp n )\sp *$. Let $\Phi$ be defined on $End(F)$ by $\Phi (\nu )=\delta \sp {-1}
\circ \varphi \circ \nu \sp {\times n} \circ \varphi \sp {-1} \circ \delta $. Since $\varphi \sp {-1}:F\sp n \to
(F\sp n )\sp * ,\; \nu \sp {\times n} :(F\sp n )\sp * \to (F\sp n )\sp * ,\;\varphi :(F\sp n )\sp * \to F\sp n $
are homomorphisms of corresponding algebras, $\Psi (\nu )= \varphi \circ \nu \sp {\times n} \circ \varphi \sp
{-1}$ is an endomorphism of $F\sp n $. Because of the condition formulated in \lemref{DeltaCondition}, this
endomorphism maps $\Delta $ into itself. Hence $\Phi (\nu )=\delta \sp {-1} \circ \Psi (\nu )\circ  \delta $ is
an endomorphism of $F$.

For two endomorphisms $\nu $ and $\mu $ of $F$ we have
$$\Phi (\nu \circ \mu )=\delta \sp {-1} \circ \varphi \circ \nu \sp {\times n} \circ \mu \sp {\times n} \circ \varphi \sp {-1} \circ \delta =$$
$$=\delta \sp {-1} \circ \varphi \circ \nu \sp {\times n} \circ \varphi \sp {-1} \circ \delta \circ  \delta \sp {-1} \circ \varphi \circ \mu \sp {\times n} \circ \varphi \sp {-1} \circ
\delta=$$
$$=\Phi (\nu ) \circ \Phi (\mu ).$$
Hence $\Phi $ is an automorphism of the semigroup $End(F)$.
\end{proof}
The permutation $\varphi$ above depends not only on automorphism $\Phi $ but on a basis of $F$. Two permutations
$\varphi$ and $\psi$ of $F\sp n$ are said to be {\it equal in effect } if they determine the same automorphism
$\Phi$ of $End(F)$. It is easy to see that this condition is equivalent to the condition that the permutation
$\varphi \circ \psi \sp {-1}$ commutes with all endomorphisms $\nu \sp {\times n}$ for $\nu \in End(F)$.

Using the description given above of every automorphism of $End(F)$ we can obtain a necessary and sufficient
condition for an automorphism $\Phi $ to be quasi-inner.

\begin{cor}\label{DirProd} An automorphism $\Phi$ of $End(F)$ is quasi-inner if and only if the corresponding permutation
$\varphi $ is equal in effect to a direct product of $n$ permutations of $F$: $\psi = s\sb 1 \times \dots \times
s\sb n$.
\end{cor}

Now we consider more useful necessary and sufficient conditions for an automorphism $\Phi $ of $End(F)$ to be
quasi-inner. Let $P$ be the subgroup of $Aut(F)$ consisting of all automorphisms $\gamma$ of $F$ such that the
restriction $\gamma$ to $X$ is a permutation of $X$, that is, $\gamma (X)=X$.
\begin{thm}\label{potinner} For any automorphism $\Phi$ of $End(F)$ the following conditions are equivalent:

1) $\Phi$ is quasi-inner;

2) $\Phi$ acts on the group $P$ in the inner way, that is, there exists an automorphism $\sigma$ of $F$ such
that $\Phi (\gamma)=\sigma \circ \gamma \circ \sigma \sp {-1}$ for every automorphism $\gamma$ from $P$;

3) there exists a basis $\tilde{X}$ of $F$ with the same number $n$ of elements such that the group $\Phi (P)$
leaves fixed $\tilde{X}$, that is $\Phi (\gamma) (\tilde{X})=\tilde{X}$ for all automorphisms $\gamma$ from $P$.
\end{thm}
\begin{proof}
$1)\Rightarrow 2)$. Let $\Phi$ is quasi-inner, that is,  for some bijection $s$ of $F$, $\Phi (\nu)=s\circ \nu
\circ s\sp {-1}$  for every endomorphism $\nu$ of $F$. Build two endomorphisms $\sigma $ and $\tau$ of $F$
setting $\sigma (x)=s(x)$ and $\tau (x)=s\sp {-1} (x)$ for all $x\in X$. Consider the composition $\mu =\Phi
(\tau )\circ \sigma $. For every $x\in X$ we have
$$\mu (x)=s\circ \tau \circ s\sp {-1}(\sigma (x))=s\circ \tau \circ s\sp {-1}(s(x))=s(\tau (x))=s(s\sp {-1} (x))=x.$$
Thus $\mu =\Phi (\tau )\circ \sigma $ is the identity map. Replacing $\Phi $ by $\Phi \sp {-1}$, we obtain that
$\Phi \sp {-1}(\sigma )\circ \tau $ is also the identity map and hence $\sigma \circ \Phi (\tau )$ is the
identity map. Conclusion: $\sigma $ is a bijection and therefore an automorphism of $F$.

Since $\sigma $ maps $X$ onto $\tilde{X}=s (X)=\sigma (X)$ the last set is also a basis of $F$. This fact will
be useful in what follows. Consider $\gamma \in P$ and calculate for every $x\in X$:
$$\sigma \sp {-1}\circ \Phi (\gamma )\circ \sigma (x)=\sigma \sp {-1}\circ s\circ \gamma \circ s\sp {-1}\circ s
(x)=\sigma \sp {-1}\circ s\circ \gamma (x)=\sigma \sp {-1}(\sigma (\gamma (x))=\gamma (x).$$

Hence $\Phi (\gamma)=\sigma \circ \gamma \circ \sigma \sp {-1}$.

$2)\Rightarrow 3)$. Now suppose that an automorphism $\sigma $ of $F$ exists, such that $\Phi (\gamma)=\sigma
\circ \gamma \circ \sigma \sp {-1}$ for every automorphism $\gamma \in P$. Set $\tilde{X}= \sigma (X)$. Clearly,
$\tilde{X}$ is a basis of $F$. For every $\gamma  \in P$ and $x\in X$, we have
$$\Phi (\gamma) (\sigma (x))=\sigma \circ \gamma \circ \sigma \sp {-1}(\sigma (x))=\sigma \circ \gamma (x),$$
hence $\Phi (\gamma)(\tilde{X})=\tilde{X}$.

$3)\Rightarrow 1)$. Let $\tilde{X} $ be a basis of $F$ satisfying the condition. Let $\sigma $ is the
automorphism which maps $X$ onto $\tilde{X}$.  Set $\Psi (\nu )=\sigma \sp {-1} \circ \Phi (\nu )\circ \sigma $
for every $\nu \in End (F)$. It is obvious that $\Psi $ is an automorphism of the monoid $End (F)$ and $\Psi
(\gamma )\in P $ for every $\gamma \in P$.

The result is that $\Phi $ is a composition of an inner automorphism and an automorphism which preserves the
group $P$. Therefore we assume that $\Phi (\gamma )\in P$ for all $\gamma  \in P$. Now consider the set $\bf C$
of all endomorphisms $\nu $ of $F$ such that $\nu \circ \gamma = \nu$ for all $\gamma \in P$. By virtue of the
condition above, $\Phi (\bf {C})=\bf {C}$. Thus $\Phi $ determines a permutation of $\bf C$. According to the
definition,  $\nu $ belongs to  $\bf C$ if and only if it has a constant value on $X$. It means that there is an
one-to-one corresponding between elements of $F$ and elements of $\bf C$, namely, for every $a\in F$ there
exists a unique endomorphism $\nu \sb a$ from $\bf C$ such that $\nu \sb a (x)=a$ for every $x\in X$.

The permutation that realizes $\Phi$ on $\bf C$ determines a permutation $s$ of $F$ in the following way:
$$s(a)=\Phi (\nu \sb a) (x) $$
for every $a\in F$.

Let $\mu$ be an endomorphism of $F$. Since $\mu \circ \nu \sb a =\nu \sb {\mu (a)}$, we obtain for every $a\in
F$
$$ \Phi (\mu )\circ \nu \sb {s(a)} = \nu \sb {s(\mu (a))}$$
and hence
$$ \Phi (\mu ) \circ s =s\circ \mu ,$$
the last gives
$$\Phi (\mu )=s\circ \mu \circ s\sp {-1}.$$
\end{proof}

\begin{defn}\label{semi-perfect} A free algebra $F$ is called {\it End-perfect} if every automorphism of $End(F)$
is inner and it is called {\it End-semi-perfect} if every automorphism of $End(F)$ is quasi-inner.
\end{defn}
It follows from the theorem above  that all  monogenic free algebras are End-semi-perfect. Now we consider free
algebras that have more than one free generators. Below, we obtain an useful sufficient condition for a free
algebra to be End-semi-perfect. A new notion we have to introduce uses $n\times n$ matrices of elements of
$n-$generated free algebra $F$. Such matrices appear as a result of the presentation of every endomorphism
$\alpha $ in a form $\alpha =\alpha \sb {(u\sb 1 , \dots ,u\sb n)}$ according to the fixed basis $x\sb 1 ,\dots
,x\sb n$. Hence we have a matrix $u\sb {(i,k)}$ for every $n$ endomorphisms $\alpha \sb i ,\; i=1,\dots ,n$
defined as follows: $u\sb {(i,k)}=\alpha \sb i (x\sb k )$.

In particular, consider $n$ special endomorphisms $\nu \sb i ,\;i=1,\dots ,n$ called {\it basis endomorphisms}
and defined as follows:
\begin{equation}\label{mainconsts}
 \nu \sb i (x\sb j )= x\sb i
\end{equation}
for all $j=1,\dots , n$. The corresponding matrix is
\[
\begin{matrix}
            x\sb 1 & x\sb 1 &\dots x\sb 1\\
            x\sb 2 & x\sb 2 &\dots x\sb 2\\
            \dots &\dots    &\dots       \\
            x\sb n & x\sb n &\dots x\sb n\\
\end{matrix}
\]

\begin{defn}\label{matrix} Let $F$ be a  $n-$generated free algebra ($n>1$). A matrix $(u\sb {i,j}), i,j=1,\dots ,n$ of its elements
is called a {\it basis matrix} if it satisfied the following conditions.

1. For every permutation $\gamma $ of indexes $1,2, \dots ,n$, there exists an automorphism $\tilde{\gamma}$ of
$F$ such that $\tilde{\gamma}(u\sb {i,j})=u\sb {\gamma (i),j}$.

2. Let $\mu \sb 1 ,\dots ,\mu \sb n $ be $n$ endomorphisms defined as follows: $\mu \sb i (x\sb j )=u\sb
{i,j},\; i,j=1,\dots , n $. Then  $\mu \sb i (u\sb {j,k})=u\sb {i,k}$ for all $i,j, k =1,\dots , n$.

3. For every $n$ endomorphisms $\alpha \sb 1 ,\dots , \alpha \sb n $ of $F$, there exists a unique endomorphism
$\alpha$ such that $\alpha (u\sb {i,j})=\alpha \sb i (u\sb {i,j})$, where $i,j=1,\dots , n$.
\end{defn}
A simple example of such a matrix was given above.

\begin{defn} An element of an algebra $F$ is said to be a {\it constant} if it is a fixed points of every endomorphism of
$F$. Denote the set of constants by $K$. An endomorphisms $\xi$ is said to be {\it constant-defined} if $\xi
(x\sb i )\in K$ for all basis elements $x\sb 1 ,\dots ,x\sb n$.
\end{defn}

It is easy to see that constant-defined endomorphisms are exactly right zeros of the semigroup $End(F)$. Denote
the set of all constant-defined endomorphisms (right zeros) by $\bf C$.  This set is invariant under every
automorphism of the semigroup $End(F)$.

\begin{prop}\label{columns} A basis matrix has the following properties.

1. It does not contain a row consisting of constants only.

2. There exists a column which does not contain two equal elements.
\end{prop}
\begin{proof}
1. Let $u\sb {ij} = w\sb 0$ be a constant. Because of the property~2 of \defnref{matrix}, every element in the
column $j$ is equal to $w\sb 0$. Thus, if a row consists of constants only, then all elements of our matrix are
constants, and what is more, all rows are equal one to another. Let $\alpha \sb 1 ,\dots , \alpha \sb n $ be
endomorphisms of $F$. Then for every endomorphism $\alpha$ we have  $\alpha (u\sb {i,j})=\alpha \sb i (u\sb
{i,j})=u\sb {i,j}$, that contraries to the property~3 of \defnref{matrix}.

2. Let all columns of our matrix contain two equal elements. Hence all elements in each column coincide. It
means that all rows coincide. By virtue of the first statement of this theorem, there is a non-constant element
$u\sb{1j}$. Thus, there exists an endomorphism $\alpha \sb 1$ such that $\alpha \sb 1 (u\sb{1j})\neq u\sb{1j}$.
Let $\alpha \sb 2 =\dots =\alpha \sb n =1\sb F$, that is, they coincide with the identity automorphism.
According to the property~3 of
\defnref{matrix}, there exists a unique endomorphism $\alpha$ such that $\alpha (u\sb {i,j})=\alpha \sb i (u\sb
{i,j})$, where $i,j=1,\dots , n$. But $\alpha (u\sb {2,j})=\alpha \sb 2 (u\sb {2,j})=u\sb {2,j}=u\sb {1,j}$ and
$\alpha (u\sb {2,j})=\alpha (u\sb {1,j})=\alpha \sb 1 (u\sb {1,j})\neq u\sb{1j}$. Contradiction.
\end{proof}
The role of basis matrixes will be explained below.
\begin{defn}
Let an automorphism $\Phi$ of $End(F)$ and a basis $x\sb 1, \dots ,x\sb n $  of $F$ be given. Consider $n$ basis
endomorphisms $\nu \sb i ,\;i=1,\dots ,n$ defined by \eqref{mainconsts}. Set $\mu \sb i =\Phi (\nu \sb i)$ and
$u\sb {i,j}=\mu \sb i (x\sb j )$. The matrix $(u\sb {i,j}), i,j=1,\dots ,n$  is called the {\it matrix of $\Phi$
with respect to the given basis}.
\end{defn}
\begin{lem}\label{automatrix}
Every matrix of an automorphism with respect to a basis is a basis matrix. Two automorphisms have the same basis
matrix with respect to the same basis if and only if one of them is product of the other and an automorphism
$\Gamma$, which leaves fixed all basis endomorphisms $\nu \sb i$ and so is quasi-inner.
\end{lem}
\begin{proof}

1. Fix a basis $x\sb 1 ,\dots , x\sb n $ of $F$. Consider $n$ basis endomorphisms $\nu \sb i ,\;i=1,\dots ,n$
defined above. It is obvious that $\nu \sb i \circ \nu \sb j =\nu \sb i $ and every permutation $\gamma $ of
indexes $1,2, \dots ,n$ produces a unique automorphism $\bar{\gamma}$ of $F$ such that  $\bar{\gamma} \circ \nu
\sb i =\nu \sb {\gamma (i)}$.

Let $\Phi$ be an automorphism of $End (F)$. Set $\mu \sb i =\Phi (\nu \sb i)$ and $u\sb {i,j}=\mu \sb i (x\sb j
)$. According to properties above, $\mu \sb i \circ \mu \sb j =\mu \sb i $ and hence $\mu \sb i (x\sb k )=\mu
\sb i (\mu \sb j (x\sb k ))$, that is, $\mu \sb i (u\sb {j,k})=u\sb {i,k}$.

Further, set $\tilde{\gamma} = \Phi (\bar{\gamma})$. We obtain $\bar{\gamma} \circ \mu \sb i =\mu \sb {\gamma
(i)}$, that is, $\tilde{\gamma}(u\sb {i,j})=u\sb {\gamma (i),j}$.

Let $\alpha \sb 1 ,\dots , \alpha \sb n $ be $n$ endomorphisms of $F$. Set $\beta \sb i =\Phi \sp {-1} (\alpha
\sb 1)$. Take $n$ elements $v\sb i =\beta \sb i (x\sb i)$. There exists a unique endomorphism $\beta $ such that
$\beta (x\sb i) =v\sb i$. Set $\alpha =\Phi (\beta )$. Since $\beta \circ \nu \sb i =\beta \sb i \circ \nu \sb
i$ under definition, we have $\alpha \circ \mu \sb i =\alpha \sb i \circ \mu \sb i$. Hence $\alpha (u\sb
{i,j})=\alpha \sb i (u\sb {i,j})$.

Thus the matrix $(u\sb {i,j}), i,j=1,\dots ,n$ is a basis matrix.

2. Let $\Phi $ and $\Psi $ be two automorphisms of $End(F)$ and their matrices with respect to a given basis
$x\sb 1 ,\dots , x\sb n $ coincide. It means that $\Phi (\nu \sb i )=\Psi (\nu \sb i )$ for all $i=1, \dots ,n$.
Consider the automorphism $\Gamma =\Psi \sp {-1} \circ \Phi $. Then $\Gamma$ lives fixed all $\nu \sb i$. Let
$\theta $ be an automorphism of $F$ from the group $P$, that is, it determines a permutation $t$ of indexes
$i=1, \dots ,n$, and hence $\theta \circ \nu \sb i =\nu \sb {t(i)}$. This implies  $\Gamma (\theta )\circ \nu
\sb i =\nu \sb {t(i)}$. The last one means that $\Gamma (\theta) \in P$. In virtue of the part~3 of
\thmref{potinner}, $\Gamma $ is quasi-inner. Thus $\Phi $ is a product $\Phi =\Psi  \circ \Gamma $ of a
quasi-inner automorphism $\Gamma $ and $\Psi$. Conversely, let the last equality take place where $\Gamma $
leaves fixed all $\nu \sb i$. It is clear that $\Phi (\nu \sb i )=\Psi (\nu \sb i )$ and their matrices
coincide.
\end{proof}
The obtained facts lead to an useful sufficient condition for an algebra to be End-semi-perfect.
\begin{thm}\label{suffcond} If the matrix of an automorphism $\Phi$ with respect to some basis of $F$ has a column
that forms a basis of $F$ then $\Phi$ is quasi-inner. Therefore if every basis matrix in $F$ has a column that
forms a basis of $F$ then $F$ is End-semi-perfect.
\end{thm}
\begin{proof} Let $\Phi$ be an automorphism and its matrix has a column $u\sb {1,j},u\sb {2,j}, \dots ,u\sb {n,j}$
which form a basis of $F$. It was mentioned in the proof of  \lemref{automatrix}, that for every $\theta \in P$
the automorphism $\Phi (\theta) $  permutes the rows of the matrix. Thus $\Phi (P)$ acts as a permutations of
the basis $u\sb {1,j},u\sb {2,j}, \dots ,u\sb {n,j}$, and $\Phi$ is quasi-inner according to \thmref{potinner}.
\end{proof}

Consider an automorphism $\Phi $ of $End(F)$ and its inverse one $\Phi \sp {-1}$. Each of them has his own basis
matrix: $(u\sb {i,j})$ and $(v\sb {i,j})$ ($ i,j=1,\dots ,n$) respectively. Define for every $m-$th column of
the first matrix an endomorphisms $\sigma \sb m$ of $F$ setting $\sigma \sb m (x\sb i )=u\sb {i,m}$. By the same
way we define for every $k-$th column of the second matrix an endomorphisms $\tau \sb k$ of $F$ setting $\tau
\sb k (x\sb i)=v\sb {i,k}$.

Denote  $\Phi (\nu \sb i)$ by $\mu \sb i$ and  $\Phi \sp {-1} (\nu \sb i)$ by $\tilde{\mu }\sb i$. Our aim is to
describe two endomorphisms $\Phi (\tau \sb k)\circ \sigma \sb m$ and $\Phi \sp {-1} (\sigma \sb m)\circ \tau \sb
k$. Start with considering the first one.
\begin{multline*}
\Phi (\tau \sb k)\circ \sigma \sb m(x\sb i)=\Phi (\tau \sb k)(u\sb {i,m})=\Phi (\tau \sb k)(\mu \sb i (x\sb m
))= \Phi (\tau \sb k\circ \nu \sb i)(x\sb m)=\\= \Phi (v\sb {i,k})(x\sb m)= \Phi (\tilde{\mu }\sb i \circ \nu
\sb k)(x\sb m)= \nu\sb i \circ\mu \sb k (x\sb m)=\nu\sb i (u\sb {k,m}).
\end{multline*}
Dually, we obtain $\Phi \sp {-1} (\sigma \sb m)\circ \tau \sb k (x\sb i )=\nu\sb i (v\sb {m,k})$.
\begin{defn}\label{pseudodiag} An endomorphism $\alpha $ of $F$ is said to be {\it pseudo-diagonal } with respect
to the fixed basis $X$ if there exists a word $w (x)$ containing only one letter $x$ such that  $\alpha
(y)=w(y)$ for every $y\in X$, that is, $\alpha $ maps $F(x)$ into itself in the same way for all basis elements
$x$.
\end{defn}
The result above can be formulated as follows.
\begin{prop}\label{sigma} The endomorphisms $\Phi (\tau \sb j)\circ \sigma \sb m$  and
$\Phi \sp {-1} (\sigma \sb m)\circ \tau \sb k$ are pseudo-diagonal or constant-defined.
\end{prop}
Now it is an interesting question if there exist numbers $k,t$ and $m$ such that endomorphisms $\Phi (\tau \sb
k)\circ \sigma \sb m$  and $\Phi \sp {-1} (\sigma \sb m)\circ \tau \sb t$ are not constant-defined. We know that
there exists element $(v\sb {i,k})$ that is not a constant. It means that $\tilde{\mu }\sb i \circ \nu \sb
k\notin {\bf C}$. Applying $\Phi $, we obtain  $\nu \sb i \circ \mu \sb k\notin {\bf C}$. It means that there
exists a number $m$ for which $\nu \sb i (u \sb {k,m})\notin K$.

For these numbers $k$ and $m$, we have $u \sb {k,m}\notin K$. It means that $\mu \sb k \circ \nu \sb m \notin
{\bf C}$. Repeating reasoning above, we find a number $t$ such that $\nu \sb k (v \sb {m,t})\notin K$. Thus we
have proved

\begin{prop}\label{notconstant} There exist such numbers $1\leq m,k,t \leq n$ for which $\nu \sb 1 (u\sb {k,m})$
and $\nu \sb 1 (v\sb {m,t})$ are not constants.
\end{prop}
\begin{thm}\label{twoMatr} Let $\Phi $ be an automorphism  of $End (F)$ such that there exist numbers
$1\leq m,k,t \leq n$ for which words  $\nu \sb 1 (u\sb {k,m})$ and $\nu \sb 1 (v\sb {m,t})$ (built from one
letter $x\sb 1$) both determine pseudo-diagonal automorphisms of $F$. Then $\Phi $ is quasi-inner.
\end{thm}
\begin{proof} Define $\sigma $ and $\tau $ as above by means of columns $(u\sb {i,m})$ and
$(v\sb {i,k})$, $i=1,\dots , n$. Under hypotheses, the pseudo-diagonal endomorphisms $\Phi (\tau \sb k)\circ
\sigma \sb m$  and $\Phi \sp {-1} (\sigma \sb m)\circ \tau \sb t$ are automorphisms. Hence $\sigma \sb m \circ
\Phi (\tau \sb t)$ also is an automorphism and the endomorphism $\sigma \sb m$ is injective and surjective,
i.e., $\sigma \sb m$ is  an automorphism. Thus $u\sb {1,m}, \dots , u\sb {n,m}$ is a basis of $F$. Now our
statement follows from \thmref{suffcond}.
\end{proof}

\section{examples}\label{examples}
It was mentioned above that all automorphisms of monogenic free algebras are quasi-inner. Now we apply the
results obtained in the previous section to prove that the most of classic free algebras are End-semi-perfect.
Some of facts presented below are known but have not been formulated in the present form and have been obtained
by a longer way.
\subsection{Semigroups and groups.}
Let $S=X\sp +$ be a free semigroup over the set $X= \{x\sb 1 ,\dots , x\sb n\}$ and $n>1$. Clearly, $\gamma
(X)=X $ for every automorphism $\gamma $ of $S$, hence  $Aut (S)=P$. Let $F$ be a free group over the set $X=
\{x\sb 1 ,\dots , x\sb n\}$ and $n>1$. It is known \cite{DyerForm} that all automorphism of $Aut (F)$ are inner,
hence all automorphisms of $End (F)$ act on the group $P$ in the inner way. According to \thmref{potinner}, we
obtain in both cases:
\begin{prop} If $F$ is a free semigroup or a free group then all automorphisms of $End (F)$ are quasi-inner.
\end{prop}
\subsection{Inverse semigroups.}
We consider inverse semigroups as algebras with two operations, a binary multiplication $\cdot$ and a unary
inversion $ \sp {-1}$ (here $a\sp {-1}$ is the inverse of an element $a$). The class of all inverse semigroups
forms a variety $\mathcal I$ defined by the identities:
$$(xy)z = x(yz),\;(xy)\sp {-1} = y\sp {-1}x\sp {-1},\;(x\sp {-1})\sp {-1} = x,$$
$$xx\sp {-1} x = x,\; x\sp {-1} xy\sp {-1} y = y\sp {-1} yx\sp {-1} x.$$

It is known that every identity $u=v$ satisfied in the variety $\mathcal I$ is a balanced one,that is, words $u$
and $v$ consists of the same letters.

Let $F=F(x\sb 1 ,\dots , x\sb n) $ be a free inverse semigroup. Denote $X = \{x\sb 1 ,\dots , x\sb n\}$, $X\sp
{-1}= \{x\sb 1 \sp {-1} ,\dots , x\sb n \sp {-1} \}$ and $X'=X\cup X\sp {-1}$. Clearly, $\theta (X')=X'$ for
every automorphism $\theta$ of $F$, and what is more, $\theta$ induces a permutation of $X'$ with the condition:
$\theta (y \sp {-1})=( \theta (y)) \sp {-1}$ for every $y\in X'$.

The left ideal $L$ corresponding to $Aut (F)$ according to introduced above Galois corresponding is the set of
all constant endomorphisms $\nu \sb e$ where $e$ is an idempotent in $F$, by the way, $L$ is the set of all left
zeros of the semigroup $End (F)$. Thus every automorphism $\Phi $ of $End (F)$ induces an automorphism of the
left ideal $L$.
\begin{prop} If $F$ is a free inverse semigroup then all automorphisms of $End (F)$ are quasi-inner.
\end{prop}
\begin{proof}
Let $(u\sb {i,j}), i,j=1,\dots ,n$ be the basis matrix corresponding to an automorphism $\Phi$ of $End (F)$.
Consider a column $k$ and suppose that two elements $u\sb {i,k}$ and $u\sb {j,k}$ contains the same letter $x\in
X$. Define endomorphisms $\alpha \sb m,\; m=1,\dots ,n$ and as follows: $\alpha \sb m (x\sb p)=x\sb m x\sb m \sp
{-1}$ for all $m,p=1,\dots ,n$. According to the property~3 in the definition of basis matrices, there exists a
unique endomorphism $\alpha $ such that $\alpha (u\sb {m,k})=\alpha \sb m (u\sb {m,k})$. Thus we have: $\alpha
(u\sb {i,k})=\alpha \sb i (u\sb {i,k}) =x\sb i x\sb i \sp {-1}$ and $\alpha (u\sb {j,k})=\alpha \sb j (u\sb
{j,k}) =x\sb j x\sb j \sp {-1}$. The first one means that $\alpha (x)$ consist from the letter $x\sb i$ only and
the second one means that $\alpha (x)$ consist from the letter $x\sb j$ only. This contradiction shows that two
elements of the same column do not contain the same letter. Conclusion: every element $u\sb {i,j}$ of the basis
matrix is built of one generator only and elements in one column are built of different letters.

$\mu \sb i =\Phi (\nu \sb i) \notin L$ because  $\nu \sb i \notin L$. It means that every row in the basis
matrix contains an element that is not an idempotent. Let $w=u\sb {1,j}$ is not an idempotent. Thus $w=
x^{-k}x^{k+m+l}x^{-l}$ or $w= x^lx^{-(k+m+l)}x^k$. In both cases $k,l\ge 0$ and $m\ge 1$. According to the
property~2 in the definition of basis matrices, $w=w^{-k}w^{k+m+l}w^{-l}$ in the first case and $w=
w^lw^{-(k+m+l)}w^k$ in the second case. Since $F$ is a combinatoric semigroup, that is, all subgroups are
trivial, and $w$ is not an idempotent, the second case does not take place and
 $m=1$ in the first case. Thus $w=w^{-k}w^{k+l+1}w^{-l}$. It is easy to see that this equality implies $k+l =0$,
that is $w=x$.

It was mentioned above that every automorphism of $F$ induces a permutation of $X'$. Thus the considered column
consists of elements of $X'$. Now apply \thmref{twoMatr} and obtain that $\Phi$ is quasi-inner.
\end{proof}
\subsection{Modules, linear spaces and Lie algebras.}
Let $F$ be a free unitary $R-$module or free Lie algebra with basis $X = \{x\sb 1 ,\dots , x\sb n\}$ over a
commutative associative ring $R$ with unit 1. Let $\Phi $ be an automorphism of $End (F)$.  Consider basis
matrices $(u\sb {i,j}), i,j=1,\dots ,n$ and $(v\sb {i,j}), i,j=1,\dots ,n$ of $\Phi $ and $\Phi \sp {-1}$
respectively. According to \propref{notconstant}, there exist numbers $k,t$ and $m$ for which $\nu \sb 1 (u\sb
{k,m})$ and $\nu \sb 1 (v\sb {m,t})$ (built from one letter $x\sb 1$) both determine pseudo-diagonal
endomorphisms. It means that $\nu \sb 1 (u\sb {k,m})=ax\sb 1$ and $\nu \sb 1 (v\sb {m,t})=bx\sb 1$ for some $a,b
\in R$ both different from zero. Since in the case $R$ is a field, these words determine automorphisms of $F$,
on the strength of \thmref{twoMatr}, we have one known and one new result (for Lie algebras):
\begin{prop} If $F$ is a finite dimensional linear space or a free finitely generated Lie algebra over a field
$K$, then all automorphisms of $End (F)$ are quasi-inner.
\end{prop}
Now assume that $R$ is a ring without zero divisors and consider the case $F$ is a free $R-$module $M$. In this
case the endomorphism $\alpha (x\sb i)=ax\sb i$ is injective. Since $\Phi (\tau \sb k )\circ \sigma \sb m
=\alpha $, $\sigma \sb m $ is injective too. Thus $u\sb {1,m},\dots ,u\sb {n,m}$ form a basis of some free
submodule $E$ of $M$. Since these modules are isomorphic, $Aut End (M)$ and $Aut End (E)$ are isomorphic too.
Repeating this process from the point $E$, we obtain a sequence $E\supseteq E\sb 1 \supseteq \dots \supseteq
E\sb k \supseteq \dots $ submodules each of which is isomorphic to $M$. If this sequence terminates on the place
$E\sb k$, then on the strength of \thmref{suffcond} all automorphisms of $End(E\sb k)$ are quasi-inner and hence
so are automorphisms of $End(M)$.
\begin{prop} If module $F$ satisfies the condition that every decreasing sequence of its $n$-generated submodules
terminates, then all automorphisms of $End (M)$ are quasi-inner.
\end{prop}

\section{derived operations}\label{DerOper}
Hereinafter, we consider only quasi-inner automorphisms of $F$. It follows from the third part of the proof of
\thmref{potinner} that we can reduce the problem to the case when an automorphism $\Phi $ preserves the group
$P$. Thus $\Phi (\mu )=s\circ \mu \circ s\sp {-1}$ for every endomorphism $\mu$, where $s(a)=\Phi (\nu \sb a)
(x) $ for every $a\in F$ because the right side of this equality does not depend on $x$. We call the permutation
defined in such way the {\it main permutation} corresponding to $\Phi$.

It was said in the first part of the mentioned proof that $s(X)$ is also a basis of $F$, and hence there is a
unique automorphism $\sigma$ coinciding with $s$ on $X$. Consequently, we have $\Phi (\gamma)=\sigma \circ
\gamma \circ \sigma \sp {-1}$ for every $\gamma \in P$. Consider the automorphism $\Psi (\mu)=\sigma \circ \mu
\circ \sigma \sp {-1}$. Since $\Psi $ coincides with $\Phi$ on $P$, $\Psi $ also preserves the group $P$.

Take the composition $\Gamma = \Psi \sp {-1} \circ \Phi$. This automorphism acts ai identity on $P$. We have
$$s \sp {\Gamma} (x)=\Gamma (\nu \sb x )(x)=\sigma \sp {-1} \circ s\circ \nu \sb x \circ s\sp {-1}\circ \sigma (x)=
\sigma \sp {-1} \circ s (x)=x.$$

This fact allows us to realize the second reducing, namely, to restrict our considerations to the case when the
main permutation preserves all basis elements.

So, $\Phi$ preserves all automorphisms from $P$ and its main permutation $s$ preserves all elements from $X$.
Let $a\sb 1 ,\dots ,a\sb n \in F$ and $\alpha \sb {a\sb 1 ,\dots ,a\sb n }$ denote the unique endomorphism which
takes $x\sb 1 $ to $a\sb 1$ ,...,$x\sb n $ to $a\sb n$ (so as in \secref{secQI}). This fact can be written in
the following way:
$$\nu \sb {a\sb i}=\alpha \sb {a\sb 1 ,\dots ,a\sb n } \circ \nu \sb {x\sb i}, \; i=1,\dots ,n .$$
Applying $\Phi$ to both sides of this equality we obtain:
$$\nu \sb {s(a\sb i)}=\Phi (\alpha \sb {a\sb 1 ,\dots ,a\sb n })\circ \nu \sb {x\sb i}, \; i=1,\dots ,n .$$
Hence
$$\Phi (\alpha \sb {a\sb 1 ,\dots ,a\sb n })=\alpha \sb {s(a\sb 1) ,\dots ,s(a\sb n) }.$$

\begin{prop} Let an endomorphism $\nu$ of $F$ satisfy the condition that for every $x\in X$ the term $\nu (x)$
can be expressed by a word that does not contain some letter $x\sb i$. Then $\Phi (\nu)$ satisfies the same
conditions.
\end{prop}
\begin{proof} Let an endomorphism $\nu$ satisfy the condition that for every $x\in X$ the term $\nu (x)$
is equal to a word which does not contain the letter $x\sb n$. This condition is equivalent to the equality
$\alpha \sb {x\sb 1 , x\sb 2 , \dots ,x\sb {n-1} ,x\sb {n-1}} \circ \nu =\nu $. Applying $\Phi$, we obtain that
$\Phi (\nu )$ satisfies the same condition.
\end{proof}
\begin{cor}Let $x\sb {i\sb 1}, \dots ,x\sb {i\sb k} \in X$. The permutation $s$ maps the free subalgebra
$F(x\sb {i\sb 1}, \dots ,x\sb {i\sb k})$ onto itself.
\end{cor}
The last result allows us to consider $\Phi$ as an automorphism of $End (F ')$ for every free subalgebra
$F'=F(x\sb {i\sb 1}, \dots ,x\sb {i\sb k})$ where $x\sb {i\sb 1}, \dots ,x\sb {i\sb k} \in X$.

Further, let $\omega$ be a $k-$ary polynomial operation in our algebra $F$, where $k\leq n$. Then the element
$s(\omega (x\sb 1,\dots ,x\sb k ))$ determines a new polynomial operation $\omega \sp *$ in $F$ as follows:
$$\omega \sp * (a\sb 1 ,\dots ,a\sb k)=\alpha \sb {a\sb 1 ,\dots ,a\sb k ,x \sb {k+1},
\dots ,x\sb n }(s(\omega (x\sb 1,\dots x\sb k ))).$$

Since for every basic operation $\omega$ we have a new (polynomial) operation $\omega \sp *$ of the same arity,
the new algebra $F\sp *$ is defined of the same type , which is a derivative algebra.
\begin{thm}\label{IsomDeriv} The main permutation $s$ is an isomorphism between $\omega $ and $\omega \sp *$.
\end{thm}
\begin{proof} Since for every $a\sb 1 ,\dots ,a\sb k \in F$
$$\omega (a\sb 1 ,\dots ,a\sb k)=\alpha \sb {a\sb 1 ,\dots ,a\sb k ,x \sb {k+1},
\dots ,x\sb n }(s(\omega (x\sb 1,\dots x\sb k ))),$$

we obtain
$$s(\omega (a\sb 1 ,\dots ,a\sb k))=\alpha \sb {s(a\sb 1) ,\dots ,s(a\sb k ),x \sb {k+1},
\dots ,x\sb n }(s(\omega (x\sb 1,\dots x\sb k )))=\omega \sp * (s(a\sb 1) ,\dots ,s(a\sb k )).$$
\end{proof}
\begin{cor} The main permutation $s$ is an isomorphism of $F$ onto $F\sp *$.
\end{cor}
The derivative algebra $F\sp *$ has the same free basis $X$. If we find all such structures, we know all
possible permutations $s$ which determine an automorphism $\Phi$. One of them may turn out to be an automorphism
of $F$. In this case $\Phi$ turns out to be inner.
\begin{defn}A permutation $\bf c$ of $F$ is called {\it central} if it commutes with all endomorphisms of $F$.
\end{defn}
\begin{thm}\label{central} An automorphism $\Phi $ is inner if and only if there exists a central permutation $\bf c$ which is an
isomorphism of $F$ onto $F\sp * $.
\end{thm}
\begin{proof}The statement can be proved in way similar to the proof of the analogous result in  \cite{PlotZhAlg}.
Let $\Phi $ is inner, then $\Phi  \sp {-1}$ is also inner, that is, it acts as follows $\Phi  \sp {-1}(\nu
)=\sigma \circ \nu \circ \sigma \sp {-1}$, where $\sigma $ is an automorphism of $F$. Set $\bf c =s\circ \sigma
$. Clearly, $\bf c $ is an isomorphism of $F$ onto $F\sp * $. Then, for every endomorphism $\nu $ of $F$ we
have: $\nu \circ {\bf c }=s\circ \Phi  \sp {-1}(\nu )\circ \sigma =s\circ \sigma \circ \nu =\bf c \circ \nu $.
Hence, $\bf c$ is a central permutation.

Conversely, let $\bf c$ is a central permutation of $F$ which is an isomorphism of $F$ onto $F\sp * $. Set
$\sigma ={\bf c} \sp {-1}\circ s$. Clearly, $\sigma$ is an automorphisms of $F$. Then,for every endomorphism
$\nu $ of $F$ we have: $\sigma \circ \nu \circ \sigma \sp {-1}={\bf c} \sp {-1}\circ s \circ \nu \circ s \sp
{-1} \circ {\bf c} ={\bf c} \sp {-1}\circ \Phi (\nu ) \circ {\bf c}=\Phi (\nu )$. Hence, $\Phi $ is inner.
\end{proof}

\section{Applications}\label{Appl}
In this section, we give some examples illustrating the suggested method. Comparing with original proofs given
by other authors shows its effectiveness. In all considered below varieties, automorphisms of $End (F)$ for
their free algebras $F$ are quasi-inner ( \secref{examples}). We use the same notation that in the previous
section.
\subsection{Semigroups}
Let $S=F(x\sb 1 , \dots ,x\sb n)$ is a free semigroup. There is only one new polynomial operation on $S$ with
respect to which $S$ is a free semigroup with the same generators: $x\sb 1 \bullet x\sb 2 =x\sb 2 x\sb 1$. Thus
the derivative structure $S\sp *$ is either the same or the dual one. Conclusion: $s$ is the identity or the
mirror bijection $\upsilon$, that assigns to every word the same word written in inverse order. The mirror
bijection determines the so called mirror automorphism of $End (S)$, described in \cite{MashSch}.
\begin{thm} Every automorphism of $End (S)$ over more than one generators is inner or a composition of an inner
automorphism and the mirror automorphism $\Upsilon$.
\end{thm}
Now consider a monogenic free semigroup $S=\{x\sp n \;\vert n\in \mathbb{N}\}$. For every $n\in \mathbb{N}\}$ we
have polynomial unary operation $\omega \sb n (a)=a\sp n$. The derived polynomial operation $\omega \sp *\sb n $
is determined by means the term $x\sp {\varphi (n)}$. Thus we have a bijection $\varphi :\mathbb{N} \to
\mathbb{N}$ and $s(a\sp n )=(s(a))\sp {\varphi (n)}$. Consequently, for every $n,m \in \mathbb{N}$, we have
$\varphi (nm)=\varphi (n)\varphi (m)$. All automorphisms of multiplicative semigroup $\mathbb{N}$ are known,
they are produced by permutations of prime numbers. Therefore we know all automorphisms of $End (S)$.
\begin{thm} Every automorphism $\Phi$ of $End (S)$ for monogenic free semigroup $S=\{x\sp n \;\vert n\in \mathbb{N}\}$
has the form $\Phi (\mu )=s\circ \mu \circ s\sp {-1}$, where $s(x\sp n) =x\sp {\varphi (n)}$ for an automorphism
$\varphi $ of the multiplicative semigroup $\mathbb{N}$.
\end{thm}
These results were first obtained in \cite{MashSch} directly.
\subsection{Inverse semigroups} Let $F$ denote a free inverse semigroup over $X$. Denote by $\upsilon $ the automorphisms
of $F$ defined as follows: $\upsilon (x)=x\sp {-1}$ for every $x\in X$. It is easy to see, that the center of
the group $Aut (F)$ consists of two elements only: the identity and $\upsilon$. Thus $\Phi (\upsilon ) =\upsilon
$. Since $\nu \sb a \circ \upsilon =\nu \sb {a\sp {-1}}$ for every $a\in F$, we obtain $s(a\sp {-1})=(s(a ))\sp
{-1} $.

Now consider the following unary operation $\omega (x)=xx\sp {-1}$. The derived operation $\omega \sp * $ is
determined by the term $s(xx\sp {-1})$. Since this term is an idempotent, it is equal to $x\sp n x\sp {-n}x\sp
{-m }x\sp m$ for some integers $n,m$, one of which does not equal to zero. Hence, for every $a\in F$, we have
$s(aa\sp {-1})=s(a)\sp n s(a)\sp {-n}s(a)\sp {-m }s(a)\sp m $. There exists an element $a$ such that $xx\sp {-1}
=s(aa\sp {-1})$. Since $xx\sp {-1}$ is a maximal idempotent in $F$, we conclude that one of integers $n$ and $m$
is equal to $0$ and the other one is equal to $1$ or to $-1$, and $a =x$. Finally we have that $s(xx\sp {-1})$
is equal to $xx\sp {-1}$ or to $x\sp {-1}x$. Result: $s(aa\sp {-1})=s(a) s(a)\sp {-1}$ for every $a\in F$ or
$s(aa\sp {-1})=s(a)\sp {-1 }s(a)$ for every $a\in F$.

Consider the case $\vert X \vert =1$, that is $F$ is a monogenic free inverse semigroup. Then $s$ acts on the
semilattice $E$ of idempotents of $F$ as an isomorphism. Let $s(aa\sp {-1})=s(a) s(a)\sp {-1}$ for all $a\in F$.
It means that $s$ leaves fixed the maximal idempotents $xx\sp {-1}$ and $x\sp {-1}x$. The conclusion is that $s$
leaves fixed the chains $\{x\sp n x\sp {-n}\}$ and $\{x\sp {-n} x\sp n \}$ where $n$ runs all positive integers.
Since $E$ is the free semilattice on the union of these ordered sets, $s$ is the identity on $E$. It is known
that $F$ is a combinatoric semigroups, and this fact implies that $s$ is the identity on $F$. If $s(aa\sp
{-1})=s(a)\sp {-1 }s(a)$ for every $a\in F$, the same reasons lead to the conclusion that $s$ is the mirror map
on $F$.

Now consider the case $n>1$. Fix two generators $x$ and $y$. Let $s(xy)=w(x,y)$, that is, the term $w(x,y)$
determines a derived binary operation. Consider the following system of three equations: $w(x,x\sp {-1})=xx\sp
{-1}$, $w (x,w (x\sp {-1}x))=x$, $w(w(xx\sp {-1}),x)=x$. The only two terms satisfying this system are  $w=xy$
and $w=yx$. \footnote {The idea of this part of the proof belongs to Gr. Mashevitzky} Thus $s(xy)=xy$ or
$s(xy)=yx$. In the first case, $s$ is the identity mapping. In the second one, $(F \sp *$ is the dual inverse
semigroup to $F$. Because the involution ${\bf c} : a\mapsto a\sp {-1}$ is an isomorphism of $F$ onto its dual
inverse semigroup $F\sp *$ and $\bf c$ is a central permutation of $F$, we conclude (\thmref{central}), that
$\Phi $ is an inner automorphism.

Thus we obtain the results first proved in \cite{MashSchZh} directly.
\begin{thm} Every automorphism of $End (F)$  is inner.
\end{thm}
\subsection{Groups}
Since the variety of all groups is a subvariety of all inverse semigroups obtained by adding the following
axiom: $x\sb 1 x\sb 1 \sp {-1}=x\sb 2 x\sb 2 \sp {-1}$, the train of thought is similar to the previous case.
Let $F$ be a free group over $X$. First of all, note that the constant endomorphism $\nu \sb e $, which assigns
to every element the group unit $e$, is a null element of $End(F)$. Hence $s(e)=e$. The same reasons that in the
case of inverse semigroups give that $s(g\sp {-1} )=(s(g))\sp {-1}$.

Now, let $\vert X \vert >1$ and $s(xy)=w(x,y)$, that is, the term $w(x,y)$ determines a derived binary operation
that gives an isomorphic free group with the same free generators. Denote the new product by $a\bullet
b=w(a,b)$. The term $w(x,y)$ has the following form: $w(x,y)=x\sp {i\sb 1}y\sp {j\sb 1}x\sp {i\sb 2}y\sp {j\sb
2}\dots x\sp {i\sb k}y\sp {j\sb k}$, where $i\sb 1 ,j\sb 1, \dots ,i\sb k ,j\sb k$ are integer numbers. The
following identities $w(x,e)=x, \;w(e,y)=y$ imply that $i\sb 1 +\dots + i\sb k =j\sb 1+\dots +j\sb k =1$, and
the identity $w(x,y)\sp {-1}=w(y\sp {-1},x\sp {-1})$ implies that $i\sb 1=j\sb k ,\dots i\sb k =j\sb 1$. Thus
$w(x,y)=x\sp {i\sb 1}y\sp {i\sb k}x\sp {i\sb 2}y\sp {i\sb {k-1} }\dots x\sp {i\sb k}y\sp {i\sb 1}$.The
consequence is the following fact $w(a,a)=a\sp 2$ for all $a\in F$. By induction, we obtain that $a\sp {\bullet
n}=a\bullet a \bullet \dots \bullet a $ ($n$ times ) is equal to $a\sp n$.

Since the source operation has a similar expression in terms of $w(x,y)$ we have an identity of the kind
$xy=x\sp {i\sb 1}\bullet y\sp {i\sb k}\bullet x\sp {i\sb 2}\bullet y\sp {i\sb {k-1} }\dots x\sp {i\sb k}\bullet
y\sp {i\sb 1}$. It can be if only if $w(x,y)$ is equal to $xy$ or to $yx$.

Since we the involution $c : g\mapsto g\sp {-1}$ plays the same role that in inverse semigroups, we conclude
\begin{thm}
Every automorphism of $End(F)$ for a free group $F$ with more than one free generators is inner.
\end{thm}
As for cyclic infinite free groups, we can repeat the train of thought for monogenic semigroups with only one
difference that instead of an automorphism $\varphi $ of the multiplicative semigroup $\mathbb{N}$ we take an
automorphism $\varphi $ of the multiplicative semigroup $\mathbb{Z}$ of all integers.

\subsection{Modules and Lie algebras.}
Let $F=F(x,y)$  be a free Lie algebra over a field $K$ considered as an algebraic system with two binary
operations, "+" and "[~]", and the set of unary operations $a\mapsto ka$ for every $k\in K$. Start with the
derivative unary operations. We have $s(kx) =\varphi (k)x$. The map $\varphi :K \to K$ is of course a bijection
preserving multiplication in $K$. Hence $\varphi (0)=0$ and $\varphi (1)=1$.

Now we assume that the field $K$ is infinite. Let   $x\bot y=s\sb F (x+y)=P(x,y)$, where $P(x,y)$ is a
polynomial in $F$. Since $x\bot y $ is homogeneous of the degree 1, the polynomial $P(x,y)$ is linear: $x\bot y
=ax+by$. It is clear that $a=b=1$ because of commutativity and of the condition $x\bot 0=x$.

Conclusion: $s $ is an additive automorphism, such that $s(kw)=\varphi (k)s(w)$, where $\varphi$ is an
automorphism of $K$, $k\in K$ and $w\in F$. Such maps are called twisted automorphisms of a module.

The next step, since $x\ast y =s\sb F ([xy])$ is a homogeneous polynomial of degree 2, $x\ast y =a[xy]$, where
$a\in K$. Consider the permutation ${\bf c}$ of $F$ defined as follows: ${\bf c}(w)=w/a$. This map is clearly
central. Hence every automorphism of  $End (F)$ is semi-inner according to definition given in
\cite{MashPlotPlot}, that is, it is induced by a twister automorphism of $F$.
\begin{thm}  Every automorphism $\Phi$ of $End (F)$ for a free Lie algebra  $F$ over an infinite field $K$
acts as follows: there exists an automorphism $\varphi $ of the field $K$ such that  $\Phi (\mu )=\sigma \circ
\mu \circ \sigma \sp {-1}$, where the bijection $\sigma$ is additive and multiplicative automorphism of the
corresponding Lie ring $F$ but $\sigma (kw)=\varphi (k)\sigma (w)$ for every $k\in K$ and $w\in F$. If $\varphi
$ is the identity automorphism then $\Phi$ is inner.
\end{thm}

We repeat the first part of the consideration above in the case of free modules $M$ over a ring $R$ satisfying
the conditions formulated in the corresponding part of the \secref{examples}, and obtain the similar result:
\begin{thm}
Every automorphism of $End(M)$ is semi-inner, that is, it is induced by a twisted automorphism of $M$.
\end{thm}

\end{document}